\documentclass[12pt]{amsart}
\usepackage{amsfonts, amssymb, amsmath, amsthm, latexsym, array}
\usepackage{fullpage,color}
\usepackage{verbatim,euscript}

\numberwithin{equation}{section}

\input {cyracc.def}
\font\tencyr=wncyr10 %scaled \magstephalf
\font\tencyi=wncyi10 %scaled \magstephalf
\font\tencysc=wncysc10 %scaled \magstephalf
\def\rus{\tencyr\cyracc}
\def\rusi{\tencyi\cyracc}
\def\rusc{\tencysc\cyracc}

\newtheorem{thm}{Theorem}[section]
\newtheorem{conj}[thm]{Conjecture}
\newtheorem{lm}[thm]{Lemma}

\newtheorem{prop}[thm]{Proposition}%[chapter]

\theoremstyle{remark}
\newtheorem{rmk}[thm]{Remark}
\newtheorem*{rem}{Remark}

\theoremstyle{definition}
\newtheorem{ex}[thm]{Example}

\newcommand {\fX}{{\mathfrak X}}

\newcommand {\slno}{{\mathfrak{sl}_{n+1}}}

\newcommand {\bi}{I}
\newcommand {\bj}{J}

\newcommand {\ca}{{\mathcal A}}

\newcommand {\co}{{\mathcal O}}

\newcommand {\BN}{{\mathbb N}}

\newcommand {\esi}{\varepsilon}

\newcommand{\ap}{\alpha}

\renewcommand{\le}{\leqslant}
\renewcommand{\ge}{\geqslant}

\newcommand{\curle}{\preccurlyeq}

\newcommand{\eus}{\EuScript}

\newcommand {\hot}{{\mathrm{ht\,}}}

\newcommand {\AN}{\mathfrak{An}}
\newcommand {\AND}{\mathfrak{An}(\Delta^+)}
\newcommand {\ANDo}{\mathfrak{An}(\Delta^+{\setminus}\,\Pi)}
\newcommand {\ANDs}{\mathfrak{An}(\Delta^+_s)}
\newcommand {\lfXr}{\langle\fX\rangle}

\newcommand {\GR}[2]{{\textrm{{\bf #1}}}_{#2}}

\newfam\eusfam%

\begin{document}
\setlength{\parskip}{3pt plus 5pt minus 0pt}
\hfill { {\color{blue}\scriptsize July 29, 2008}}
\vskip1ex

\title[On orbits of antichains of positive roots]
{On orbits of antichains of positive roots}
\author[D.\,Panyushev]{Dmitri I.~Panyushev}
\address[]{Institute for Information Transmission Problems, B. Karetnyi per. 19, Moscow 101447\hfil\break\indent
Independent University of Moscow,
Bol'shoi Vlasevskii per. 11, 119002 Moscow, \ Russia}
\email{panyush@mccme.ru}
\thanks{Supported in part by  R.F.B.R. grant
06--01--72550.}
\begin{abstract}
For any finite poset $\eus P$, there is a natural operator, $\fX=\fX_{\eus P}$,
acting on the set of antichains of $\eus P$.
We discuss conjectural properties of $\fX$
for some graded posets associated with irreducible root systems. 
In particular, if $\Delta^+$ is the set of positive roots and  $\Pi$ is the set of simple
roots in $\Delta^+$, then we consider the cases $\eus P=\Delta^+$ and
$\eus P=\Delta^+\setminus\Pi$. For the root system of type $\GR{A}{n}$, we 
consider an $\fX$-invariant integer-valued function on the set of antichains of $\Delta^+$
and establish some properties of it.
\end{abstract}
\maketitle

\section{Introduction}  

\noindent 
Let $(\eus P,  \curle)$ be an arbitrary finite poset. 
For any  $\eus S\subset \eus P$, let $\eus S_{min}$ and $\eus S_{max}$ denote
the set of minimal and maximal elements of $\eus S$, respectively.
An {\it antichain\/} in $\eus P$ is a 
subset of mutually incomparable elements. In other words, $\Gamma$ is an antichain if
and only if $\Gamma=\Gamma_{min}$ (or $\Gamma=\Gamma_{max}$).
Write $\AN(\eus P)$ for the set of all antichains in $\eus P$. 
An {\it upper ideal\/} (or {\it filter}) is a subset $\eus I\subset\eus P$ such that if $\gamma\in \eus I$ and 
$\gamma\curle \beta$, then $\beta\in \eus I$. If $\Gamma\in \AN(\eus P)$, 
then $\eus I(\Gamma)$ denotes the { upper ideal\/} of $\eus P$ generated by $\Gamma$. 
That is, 
\[
   \eus I(\Gamma)=\{ \esi\in \eus P\mid  \exists \gamma\in \Gamma \text{ such that } 
   \gamma\curle\esi \} \ .
\]
For instance, $\Gamma=\varnothing$ is an antichain and $\eus I(\varnothing)$ is 
the empty upper ideal.
Conversely, if $\eus I$ is an upper ideal of $\eus P$,
% and $\Gamma(\eus I)$ is the set of minimal elements of $\eus I$, 
then $\eus I_{min}\in\AN(\eus P)$.
This yields a natural bijection between the upper ideals and antichains of $\eus P$.
Letting $\Gamma'\lessdot \Gamma$ if $\eus I(\Gamma')\subset \eus I(\Gamma)$, we make
$\AN(\eus P)$ a poset. 

For $\Gamma\in \AN(\eus P)$,  we set $\fX(\Gamma)=
(\eus P\setminus \eus I(\Gamma))_{max}$. This
defines the map $\fX=\fX_{\eus P}: \AN(\eus P)\to \AN(\eus P)$.
Clearly,   $\fX$ is one-to-one, i.e., it is a
permutation of the finite set $\AN(\eus P)$. We say that $\fX$ is the {\it reverse operator}
for $\eus P$. If $\#\AN(\eus P)=m$, then $\fX$ is an element of the symmetric group $\Sigma_m$.
Let $\lfXr$ denote the cyclic subgroup of $\Sigma_m$ 
generated by $\fX$. The {\it order\/} of $\fX$, $\text{ord}(\fX)$, is the order of  the group $\lfXr$.
As the definition of $\fX$ is quite natural, one can expect that properties of $\lfXr$-orbits
in $\AN(\eus P)$ are closely related to other properties of $\eus P$. 
One of the problems is to determine the cyclic structure of $\fX$, i.e., possible cardinalities 
of  $\lfXr$-orbits  in $\AN(\eus P)$.
In particular, one can ask about a connection between properties of 
$\eus P$ and $\text{ord}(\fX)$.  For simplicity,
we will speak about $\fX$-orbits in what follows.
If $\eus P$ is a Boolean lattice, then  $\fX$-orbits has been studied
under the name "loops of clutters", see \cite{clutters}. Some conjectures 
stated in  \cite{clutters} for that special situation are proved  in \cite{flaas1,flaas2} for an 
arbitrary graded poset $\eus P$.

We say that $\eus P$ is {\it graded (of level\/} $r$) if there is a function 
$d:\eus P\to \{1,2,\dots,r\}$ such that both $d^{-1}(1)$ and $d^{-1}(r)$ are non-empty, and
$d(y)=d(x)+1$ whenever  $y$ covers $x$. Then $d^{-1}(1)\subset \eus P_{min}$ and 
$d^{-1}(r)\subset \eus P_{max}$.

\begin{lm}   \label{stand_orb}
Suppose $\eus P$ is graded of level $r$,  $d^{-1}(1)= \eus P_{min}$ and 
$d^{-1}(r)= \eus P_{max}$.
%has a unique maximal element.
Then $\fX$ has an orbit of cardinality $r+1$. 
\end{lm}\begin{proof}
Clearly, $\eus P(i):=d^{-1}(i)$ is an antichain for any $i$. From our hypotheses, 
it follows that $\fX(\eus P(i))=\eus P(i{-}1)$ for $i=2,\dots,r$, 
$\fX(\eus P(1))=\varnothing$, and $\fX(\varnothing)=\eus P(r)$.
Thus,
$\{\varnothing, \eus P(r),\dots,\eus P(1)\}$ is an $\fX$-orbit.
\end{proof}

Such an orbit of $\fX$ is said to be {\it standard}.

\noindent
The goal of this note is to present several observations and conjectures 
on orbits of reverse operators for some graded posets associated with a root system
$\Delta$.
In Section~2, we discuss conjectural properties of reverse operators for
$\Delta^+$, $\Delta^+{\setminus}\,\Pi$, and $\Delta^+_s$ (see notation below). 
Roughly speaking, all our conjectures are verified up to rank 5. In particular, our calculations
for $\GR{F}{4}$ are presented in Appendix.
In Section~3, we work with the root system of type $\GR{A}{n}$.
In this particular case, we 

(1)  describe an $\fX$-invariant function $\eus Y:\AND\to \BN$ (this is due to O.\,Yakimova);

(2) prove that $\fX$ satisfies the relation $\fX^{-1}=\ast\circ\fX\circ\ast$, where 
$\ast:\AND\to \AND$ is the involutory mapping (duality) constructed in \cite{dual}.
In other words, for any $\Gamma\in\AND$, one has $\fX^{-1}(\Gamma^*)=\fX(\Gamma)^*$;

(3) show that $\eus Y(\Gamma)=\eus Y(\Gamma^*)$ for any $\Gamma\in\AND$.

\vskip1ex\noindent
This is an expanded version of my talk at the workshop
``$B$-stable ideals and nilpotent orbits" (Rome, October 2007).

%%%%%%%%%%%%%
\section{Reverse operators for posets associated with root systems}  \label{sect:Delta}
%%%%%%%%%%%%%

\noindent 
Let $\Delta$ be a reduced irreducible root system
in an $n$-dimensional real vector space $V$ and $W\subset GL(V)$ the corresponding
Weyl group.
Choose a system of positive roots $\Delta^+$ with the corresponding
subset of simple roots $\Pi=\{\ap_1,\ldots,\ap_n\}$. 
The {\it root order\/} in $\Delta^+$ is given  by letting $x\curle y$ if
$y-x$ is a non-negative integral combination of positive roots. 
%%This defines a partial order in $V$. 
In particular, $y$ covers $x$ if $y-x$ is a simple root.  The highest root in $\Delta^+$ 
is denoted by $\theta$. It is the unique maximal element of $(\Delta^+, \curle)$. 
If $\Delta$ has two root lengths, then $\theta_s$ is the dominant (highest) short root.
Let $w_0\in W$ be the longest element, i.e., the unique element that takes $\Delta^+$ to
$-\Delta^+$. If $\gamma=\sum_{i=1}^n a_i \ap_i\in\Delta^+$, then $\hot(\gamma):=\sum a_i$
is the {\it height\/} of $\gamma$. For $I\subset \Pi$,  $\Delta(I)$ is the root subsystem of 
$\Delta$ generated by $I$. If $\GR{X}{n}$ is one of the Cartan types, then
$\Delta(\GR{X}{n})$ denotes the root system of type $\GR{X}{n}$.

\subsection{Orbits in $\Delta^+$}
In this subsection, we consider antichains in $\Delta^+$ and the 
reverse operator $\fX=\fX_{\Delta^+}:\AND \to \AND$.
 
Let  $h=h(\Delta)$ be the Coxeter number and $e_1,\dots,e_n$ the exponents of 
$\Delta$.  It is known \cite{cp2} that $\#(\AND)=\displaystyle
\prod_{i=1}^n \frac{h+e_i+1}{e_i+1}$.
The function
$\ap\mapsto\hot(\ap)$ turns $\Delta^+$ into a graded poset of level $h{-}1$.
Set $\Delta(i)=\{\ap\in\Delta^+\mid \hot(\ap)=i\}$
and $\Delta({\ge}i)=\{\ap\in\Delta^+\mid \hot(\ap)\ge i\}$. Then
$\Delta(1)=\Pi=\Delta^+_{min}$ and $\Delta(h{-}1)=\{\theta\}=\Delta^+_{max}$.

Let us point out two specific orbits of $\fX$:

1) By Lemma~\ref{stand_orb}, there is an orbit of cardinality $h$. Namely,
$\{ \varnothing, \Delta(h{-}1),\dots,\Delta(2),\Delta(1)\}$ is the standard 
$\fX$-orbit in $\AND$.

2)  There is an $\fX$-orbit of order $2$. Let $\ca\subset\Pi$ be a set of mutually
orthogonal roots such that
$\Pi\setminus \ca$ also has that property. (The partition $\{\ca,\Pi\setminus\ca\}$
is uniquely determined, since the Dynkin diagram of $\Delta$ is a tree.) Then 
$\fX(\ca)=\Pi\setminus\ca$ and $\fX(\Pi\setminus\ca)=\ca$.

If $\Delta$ is of rank 2, then these two orbits exhaust $\AND$.

\begin{conj}   \label{conjD}
{\ } \phantom{aaa} 
\begin{itemize}
\item[\sf (i)] \ 
If $w_0=-1$, then ${\rm ord}(\fX)=h$;
\item[\sf (ii)] \  If  $w_0\ne -1$, then $\fX^h$ is the involution 
of\/ $\AND$ induced by $-w_0$ and ${\rm ord}(\fX)=2h$;
\item[\sf (iii)] \ Let  $\co$ be an arbitrary $\fX$-orbit in $\AND$.
Then $\displaystyle \frac{1}{\#\co} \sum_{\Gamma\in \co}  \# \Gamma=
\frac{\# \Delta^+}{h}=\frac{n}{2}$.
\end{itemize}
\end{conj}

\noindent Recall that $w_0\ne -1$ if and only if $\Delta$ is of type 
$\GR{A}{n}$ ($n\ge 2$), $\GR{D}{2n+1}$, $\GR{E}{6}$. Furthermore, the posets
$\Delta^+$ are isomorphic for $\GR{B}{n}$ and $\GR{C}{n}$ \cite[Lemma\,2.2]{coveri}.
Conjecture~\ref{conjD} has been verified for  $\GR{A}{n}$ ($n\le 5$), $\GR{C}{n}$ ($n\le 4$),
$\GR{D}{4}$, $\GR{F}{4}$. 
It is easily seen that  $\#\Gamma$ equals the number of elements of $\AND$ covered by 
$\Gamma$. For, $\Gamma$ covers $\Gamma'$ with respect to the order `$\lessdot$'
described in the Introduction if and only if
$\Gamma'=(\eus I(\Gamma)\setminus\{\gamma_i\})_{min}$ for some
$\gamma_i\in\Gamma$.
Hence $\displaystyle \sum_{\Gamma\in \AND}\# \Gamma$ equals the total number of
edges in the Hasse diagram of $(\AND, \lessdot)$. Therefore  it follows from 
\cite[Cor.\,3.4]{jac06} that 
\[
  \sum_{\Gamma\in \AND} \frac{\# \Gamma}{\#\AND} =\frac{\# \Delta^+}{h} \ .
\]
Thus, part (iii) can be regarded as a refinement of the last equality.

\begin{ex}
We use the standard notation for roots in $\Delta^+(\GR{A}{n})$; e.g.,
$\ap_i=\esi_i-\esi_{i+1}$, $i=1,2,\dots,n$, and $\theta=\esi_1-\esi_{n+1}$.
If $\Gamma=\{\ap_1\}$ for $\GR{A}{n}$ and $n\ge 3$, then 
\[
\fX^k(\{\ap_1\})=\{ \gamma\in \Delta(\ap_1,\dots,\ap_{n-1})\mid \hot(\gamma)=n+1-k \}
\sqcup \{\ap_{k+1}+\ldots +\ap_n\}, \ 1\le k\le n.
\]
In particular, $\fX^n(\{\ap_1\})=\{\ap_1,\dots,\ap_{n-1}\}$ and hence 
$\fX^{n+1}(\{\ap_1\})=\{\ap_n\}$. Therefore the $\fX$-orbit of $\{\ap_1\}$ is of 
cardinality $2h=2n+2$. The ratio 
$\displaystyle \frac{1}{\#\co} \sum_{\Gamma\in \co}  \# \Gamma$ equals 
$n/2$ for this orbit, as required.
\end{ex}

%Appearance of the Coxeter number $h$ suggests that $\fX$ could somehow be related  
%to a Coxeter element of the Weyl group  $W$. 

It is an interesting problem to construct ``invariants'' of $\fX$, i.e., functions on
$\AND$ that are constant on the $\fX$-orbits. Ideally, one could ask for a family of invariants 
that separates the orbits. Our achievement in this direction is rather modest.
We know only one invariant in the case of type $\GR{A}{n}$, see Section~\ref{sect:An}.

\subsection{Orbits in $\Delta^+{\setminus}\,\Pi$}
We regard $\Delta^+{\setminus}\,\Pi=\Delta({\ge}2)$ as subposet of $\Delta^+$. 
The theory of antichains (upper ideals) in $\Delta^+{\setminus}\, \Pi$ is quite similar to that
for $\Delta^+$ \cite{eric}. In particular,
$\#(\ANDo)=\displaystyle \prod_{i=1}^n \frac{h+e_i-1}{e_i+1}$.
Let $\fX_0: \ANDo \to \ANDo$ be  the reverse operator  for 
$\Delta^+{\setminus}\, \Pi$.
The function $\ap\mapsto (\hot\ap){-}1$ turns $\Delta^+{\setminus}\,\Pi$ into a graded poset of
level $h{-}2$. It follows that  $\fX_0$  has the standard orbit of cardinality $h{-}1$.
As the simple roots are removed, the corresponding orbit of order 2 also vanishes from
$\ANDo$.

\begin{conj}   \label{conjD0}
{\ } \phantom{aaa} 
\begin{itemize}
\item[\sf (i)] \ 
If $w_0=-1$, then ${\rm ord}(\fX_0)=h-1$;
\item[\sf (ii)] \  If  $w_0\ne -1$, then $\fX_0^{h-1}$ is the involution 
of\/ $\ANDo$ induced by $-w_0$ and ${\rm ord}(\fX_0)=2h{-}2$;
\item[\sf (iii)] \ For any   $\fX_0$-orbit $\co\subset\ANDo$, we have
$\displaystyle \frac{1}{\#\co}\sum_{\Gamma\in \co}  \# \Gamma=
\frac{\# (\Delta^+\setminus \Pi)}{h-1}=\frac{n}{2}{\cdot}\frac{h{-}2}{h{-}1}$\, .
\end{itemize}
\end{conj}

\noindent
Here are empirical evidences supporting the conjecture.
The poset $\Delta^+{\setminus}\, \Pi$ for $\GR{A}{n+1}$ is isomorphic to 
$\Delta^+$ for $\GR{A}{n}$. Therefore Conjecture~\ref{conjD0} holds for 
$\GR{A}{n}$ ($n\le 6$).
It has also been verified for $\GR{C}{n}$ ($n\le 5$), $\GR{D}{n}$ ($n\le 5$), 
and  $\GR{F}{4}$. 
Again, $\displaystyle \sum_{\Gamma\in \ANDo}\# \Gamma$ equals the number of
edges on the Hasse diagram of $\ANDo$, and it was verified in  \cite[Sect.\,3]{coveri}
that 
\[
 \frac{1}{\#\ANDo}\sum_{\Gamma\in \ANDo} \#\Gamma =\frac{\# (\Delta^+\setminus\Pi)}{h-1}\ .
\]
Hence part (iii) can be regarded as a refinement of the last equality.

If $w_0=-1$ and $h{-}1$ is prime, then  Conjecture~\ref{conjD0} predicts  that
all $\fX_0$-orbits have the same cardinality. This is really the case for $\GR{F}{4}$,
$\GR{C}{3}$, and $\GR{C}{4}$. Actually, this seems to be true for any $\GR{C}{n}$, see
Conjecture~\ref{C-short}.

\begin{rem}
One might have thought that posets $\Delta({\ge}j)$ enjoy similar good properties for any
$j$. However,
this is not the case. For $\GR{F}{4}$ and $\Delta({\ge}3)$, the reverse operator 
has orbits of cardinality 10 and 8. Hence its order equals 40, while $h-2=10$. 
Furthermore, the mean value of the
size of antichains along the orbits is not constant.
\end{rem}

\subsection{Orbits in $\Delta^+_s$}   \label{sect:delta_s}
%\noindent
Suppose $\Delta$ has two root lengths.
Then $\Delta^+_s$ denotes the set of short positive roots in $\Delta^+$.
We regard $\Delta^+_s$ as subposet of $\Delta^+$. Then $\theta_s$ is the unique maximal
element of $\Delta^+_s$ and $(\Delta^+_s)_{min}=\Pi\cap \Delta^+_s=:\Pi_s$.
General results on antichains in $\Delta^+_s$ are obtained in \cite[Sect.\,5]{losh}.
Suppose $m=\#\Pi_s$ and the exponents $\{e_i\}$ are increasingly ordered.
Then $\#(\ANDs)=\displaystyle \prod_{i=1}^m \frac{h+e_i+1}{e_i+1}$.
Let  
$\fX_s:\ANDs\to \ANDs$ be the reverse operator for $\ANDs$.
Let $h^*(\Delta)$ denote  the {\it dual Coxeter number\/} of $\Delta$.
Recall that $h^*(\Delta^\vee)-1=\hot(\theta_s)$, where $\Delta^\vee=
\{ \frac{2\ap}{(\ap,\ap)} \mid \ap\in\Delta\}$ is the dual root system.
The function $\hot(\ )$ turns $\Delta^+_s$
into a graded poset of level $h^*(\Delta^\vee)-1$. It follows that 
$\fX_s$  has the standard orbit of cardinality $h^*(\Delta^\vee)$.

\begin{conj}   \label{conjDs}
{\ } \phantom{aaa} 
\begin{itemize}
\item[\sf (i)] \ 
${\rm ord}(\fX_s)=h^*(\Delta^\vee)$;
\item[\sf (ii)] \ Let  $\co$ be an arbitrary $\fX_s$-orbit in $\ANDs$.
Then $\displaystyle \frac{1}{\#\co}\sum_{\Gamma\in \co} \# \Gamma =
\frac{\# (\Delta^+_s)}{h^*(\Delta^\vee)}$.
\end{itemize}
\end{conj}

\noindent
The conjecture is easily verified for $\GR{B}{n}$, $\GR{F}{4}$, and $\GR{G}{2}$, where 
the number of $\fX_s$-orbits equals 1,\,3, and 1, respectively. We have also 
verified it for $\GR{C}{n}$ with $n\le 5$. 

For $\GR{C}{n}$, the posets $\Delta^+{\setminus}\,\Pi$ and $\Delta^+_s$ (hence
$\ANDo$ and  $\ANDs$) are isomorphic.
We also have a more precise conjecture in this case.

\begin{conj}  \label{C-short}
For\/ $\Delta^+_s(\GR{C}{n})$, every $\fX_s$-orbit is of cardinality $2n-1=h^*(\GR{B}{n})$.
Each $\fX_s$-orbit contains a unique antichain lying in 
$\Delta^+(\ap_1,\dots,\ap_{n-2})\simeq \Delta^+(\GR{A}{n-2})$.
\end{conj}

Since $\#(\ANDs)=\genfrac{(}{)}{0pt}{}{2n-1}{n}$ for $\GR{C}{n}$ 
\cite[Theorem\,5.5]{losh}, 
Conjecture~\ref{C-short} would imply
that the number of $\fX_s$-orbits equals $\frac{1}{2n-1}\genfrac{(}{)}{0pt}{}{2n-1}{n}$,
the $(n{-}1)$-th Catalan number. Note that this conjecture also provides a canonical representative in each $\fX_s$-orbit in $\AN(\Delta^+_s(\GR{C}{n}))$.

%%%%%%%%%%%%%%
\subsection{Orbits in $\Delta^+_s\setminus \Pi_s$}   \label{sect:delta_s-minus}
\phantom{ } \ \  
We regard $\Delta^+_s\setminus \Pi_s$ as subposet of $\Delta^+_s$. 
For the reverse operator 
$\fX_{s,0}: \AN(\Delta^+_s\setminus \Pi_s)\to \AN(\Delta^+_s\setminus \Pi_s)$, one can
state a similar conjecture, where $h^*(\Delta^\vee)$ is replaced with 
$h^*(\Delta^\vee)-1$. However, this does not make much sense. The case of
$\GR{B}{n}$ and $\GR{G}{2}$ is trivial.
For $\GR{C}{n}$, the poset
$\Delta^+_s\setminus \Pi_s$ is isomorphic to $\Delta^+(\GR{C}{n-1})$. 
Hence this case is covered by previous conjectures. The only new phenomenon occurs
for $\GR{F}{4}$, where everything is easily verified. 
Here $\# \AN(\Delta^+_s\setminus \Pi_s)=16$ and 
$\fX_{s,0}$ has two orbits, both of cardinality $8=h^*(\GR{F}{4})-1$.

\begin{ex} A slight modification of a poset can drastically change 
properties of reverse operators.
Consider two graded posets of level 3, with Hasse diagrams
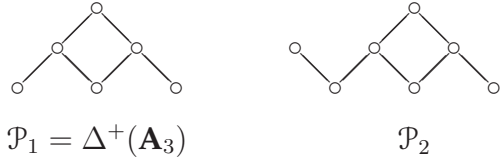
\begin{center}
\begin{picture}(210,75)(30,-285)
\multiput(60,-250)(30,0){3}{\circle{4}}
\multiput(75,-235)(30,0){2}{\circle{4}}
\put(90,-220){\circle{4}}
\multiput(77,-237)(30,0){2}{\line(1,-1){11}}
\multiput(73,-237)(30,0){2}{\line(-1,-1){11}}
\put(92,-222){\line(1,-1){11}}
\put(88,-222){\line(-1,-1){11}}
\put(85,-275){\makebox(10,10){$\eus P_1=\Delta^+(\GR{A}{3})$}}

\multiput(180,-250)(30,0){3}{\circle{4}}
\multiput(165,-235)(30,0){3}{\circle{4}}
\put(210,-220){\circle{4}}
\multiput(167,-237)(30,0){3}{\line(1,-1){11}}
\multiput(193,-237)(30,0){2}{\line(-1,-1){11}}
\put(212,-222){\line(1,-1){11}}
\put(208,-222){\line(-1,-1){11}}
\put(205,-275){\makebox(10,10){$\eus P_2$}}
\end{picture}
\end{center}
The reverse operator for $\eus P_1$  has three orbits of cardinality 8,4, and 2
(and the properties stated in Conjecture~\ref{conjD}). 
For $\eus P_2$, there are two orbits of cardinality 16 and 7.
Thus, $\text{ord}(\fX_1)=8$, while $\text{ord}(\fX_2)=16{\cdot}7$. Furthermore, the mean 
values of the size of antichains for two $\fX_2$-orbits are different.
\end{ex}

%%%%%%%%%%%%%
\section{Results for $\Delta^+(\GR{A}{n})$}   \label{sect:An}
%%%%%%%%%%%%%

\noindent In this section, 
$\Delta=\Delta(\GR{A}{n})=\Delta(\slno)$.

\subsection{The OY--invariant}  Here we describe an $\fX$-invariant function  
$\eus Y:\AND\to \BN$, which is found by Oksana Yakimova.

Let $\Gamma=\{\gamma_1,\dots,\gamma_k\}$ be an arbitrary antichain in $\Delta^+$
and $\eus I=\eus I(\Gamma)$ the corresponding upper ideal, so that $\Gamma=\eus I_{min}$.
To each $\gamma_s$, we attach certain integer as follows.
Clearly, $\eus I\setminus\{\gamma_s\}$ is again an upper ideal. 
%denote the subset of its minimal elements.
Set 
\[
     r_{\Gamma}(\gamma_s):=\#(\eus I\setminus\{\gamma_s\})_{min}-\# \eus I_{min} +1 \ .
\]
For $\slno$, the difference between the numbers of minimal elements 
of $\eus I$ and $\eus I\setminus\{\gamma_s\}$ always 
belongs to $\{-1,0,1\}$. Therefore $r_{\Gamma}(\gamma_s)\in\{0,1,2\}$. The
OY-number of $\Gamma$ is defined by
\begin{equation}   \label{def-oy}
   \eus Y(\Gamma):=\sum_{s=1}^k r_{\Gamma}(\gamma_s) .
\end{equation}
This definition only applies to non-empty $\Gamma$, and
we specially set $\eus Y(\varnothing)=0$.

\begin{ex}  \label{ex:valuesY}
{\sf a)} For $\Gamma=\Pi=\{\ap_1,\dots,\ap_{n}\}$, we have $\eus Y(\Pi)=0$. More generally, 
the same is true for $\Gamma=\Delta(i)$. \ 
{\sf b)} For $\Gamma=\{\ap_1,\ap_3,\dots \}$ (all simple roots with odd numbers)
or $\Gamma=\{\ap_2,\ap_4,\dots \}$ (all simple roots with even numbers), we have
$\eus Y(\Gamma)=n-1$.
\end{ex}
\begin{thm}[O.\,Yakimova]           \label{thm:osya}
The OY-number is $\fX$-invariant, i.e., $\eus Y(\Gamma)=\eus Y(\fX(\Gamma))$ for all\/
$\Gamma\in\AND$.
\end{thm}\begin{proof}
Let us begin with an equivalent definition of $\eus Y(\Gamma)$, which is better for  
the proof. Recall that $\Delta^+(\GR{A}{n})=\{ \esi_i-\esi_{j+1} \mid
1\le i\le j\le n\}$.
The positive root $\esi_i-\esi_{j+1}=\ap_i+\ldots +\ap_j$ will be denoted by $(i,j)$.
Suppose $\gamma_s=(i_s,j_s)$. Without loss of generality, we may assume that the 
$i$-components of all roots in $\Gamma$ form an increasing sequence. 
Then the fact that $\Gamma=\{(i_1,j_1),\dots,(i_k,j_k)\}$ is an antichain is equivalent to that 
$1\le i_1< \ldots < i_k$, \ $j_1< \ldots < j_k\le n$, and $i_s \le j_s$ for each $s$.
Obviously, 
$\Gamma\setminus\{\gamma_s\}\subset (\eus I\setminus\{\gamma_s\})_{min}$.
%$\Gamma\setminus \{\gamma_s\}$ is contained in the set of  minimal elements
%of  the upper ideal $\eus I\setminus \{\gamma_s\}$. 
Furthermore, if $i_s{-}i_{s-1}\ge 2$, then
$(i_s{-}1,j_s)\in (\eus I\setminus\{\gamma_s\})_{min}$;  
and if $j_{s+1}{-}j_s\ge 2$, then
$(i_s,j_s{+}1)\in (\eus I\setminus\{\gamma_s\})_{min}$ as well.
This observation shows that
 $r_\Gamma(\gamma_s)=\chi(i_s{-}i_{s-1})+\chi(j_{s+1}{-}j_s)$, where 
 $i_0:=0$, $j_{k+1}:=n+1$,
and  the function $\chi$ on $\{1,2,\dots\}$ is defined  by 
\[  \chi(a)=\begin{cases}
1, &  a\ge 2 \\ 0, & a=1 \end{cases}.
\]
Hence
\begin{equation}  \label{conv_oy}
  \eus Y(\Gamma)=\sum_{s=1}^k\chi(i_s{-}i_{s-1}) +
  \sum_{s=1}^k \chi(j_{s+1}{-}j_s) .
\end{equation}
We say that the difference $b-a$ is {\it essential\/} if $b-a\ge 2$.
Thus, $\eus Y(\Gamma)$ counts the total number of consecutive 
essential differences  in the sequences 
$(0,i_1,\dots,i_k)$ and $(j_1,\dots,j_k,n+1)$. For this reason,
we will think of $\Gamma$ as two-row array:
\begin{equation}   \label{eq:gamma}
    \Gamma=\left(\begin{array}{ccccc}0 & i_1 & \dots & i_k & \\ 
    & j_1 & \dots & j_k & n+1\end{array}\right),
\end{equation}
where each 2-element column represents a positive root.
\begin{figure}[htb]
\setlength{\unitlength}{0.04in}
\begin{center}
\begin{picture}(95,56)(-5,35)

\put(5,90){\line(1,0){75}}    %%%  horizontal
\put(80,35){\line(0,1){55}}          
\linethickness{.3mm}
{\color{magenta}\put(20,70){\line(1,0){15}}                   
\put(35,55){\line(1,0){20}}                   
\put(55,45){\line(1,0){24}}
                   
\put(20,70){\line(0,1){20}}       %% vertical
\put(35,55){\line(0,1){15}}                   
\put(55,45){\line(0,1){10}} }                   

\thinlines
\put(20,74){\line(1,0){4}}  
\put(35,59){\line(1,0){4}}  
\put(24,70){\line(0,1){4}}                   
\put(39,55){\line(0,1){4}}      
\put(32,82){\vector(-1,-1){10}}     \put(33,81){\footnotesize $(i_1,j_1)$}     %% generators
\put(47,67){\vector(-1,-1){10}}     \put(48,66){\footnotesize $(i_2,j_2)$}       

{\color{cyan}
\put(16,86){\line(1,0){4}}  
\put(31,66){\line(1,0){4}}  
\put(75,41){\line(1,0){4}}  
\put(16,86){\line(0,1){4}}                   
\put(31,66){\line(0,1){4}}  
\put(75,41){\line(0,1){4}} }    

\put(8,78){\vector(1,1){10}}     \put(-4,75){\footnotesize $(1,j_1{-}1)$}     %% generators
\put(23,58){\vector(1,1){10}}     \put(8,55){\footnotesize $(i_1{+}1,j_2{-}1)$}       
\put(67,33){\vector(1,1){10}}     \put(58,30){\footnotesize $(i_k{+}1,n)$}       

\end{picture}
\end{center}
\caption{Antichains $\Gamma$ and $\fX(\Gamma)$ for $\slno$}   \label{pikcha_sl}
\end{figure}

\noindent
Let us describe the operator $\fX$ using this notation. The first step is to replace 
$\Gamma$ in Eq.~\eqref{eq:gamma} with
\begin{equation}   \label{eq:new_gamma}
    \tilde\Gamma=\left(\begin{array}{ccccccc}
   0   &   1         & i_1{+}1 & \dots &  i_{k-1}{+}1 & i_k{+}1 &  \\ 
        &  j_1{-}1 & j_2{-}1  & \dots &  j_k{-}1        & n   & n{+}1 
\end{array}\right).
\end{equation}
It may happen, however, that some 2-element  
columns of $\tilde\Gamma$ are ``bad'', i.e., they do not represent positive 
roots; e.g., if $j_1=1$ or $i_{s-1}{+}1 > j_s{-}1$. The second step is to remove all 
bad columns. The remaining array is exactly $\fX(\Gamma)$, cf. Figure~\ref{pikcha_sl}.

Thus, our task is to check that such a procedure 
does not change the total number of essential differences.

(a)  \ If $\fX(\Gamma)=\tilde\Gamma$, then the essential differences themselves
for  $\Gamma$ and $\fX(\Gamma)$ are the same.

(b) \ Let us realise what happens with essential differences if $\tilde\Gamma$
contains bad columns. Assume the column 
$\varkappa_s=\left(\!\begin{array}{c} i_{s-1}{+}1 \\ j_s{-}1
\end{array}\!\right)$ is bad for $2\le s\le k-1$. 
It is easily seen that in this case $i_{s-1}{+}1= j_s$
and $\gamma_{s-1}$, $\gamma_{s}$ are adjacent simple roots.
If both the surrounding columns for $\varkappa_s$ are good and, say, 
$\gamma_s=\ap_t=(t,t)$, 
then the array $\tilde\Gamma$ contains a fragment of the form
\[
  \left(\begin{array}{ccccc} \dots & x    & t+1 & t+2 & \dots \\
                                          \dots & t-1 & t      & y    & \dots
\end{array}\right), 
\]  
where $x\le t-1$ and $y\ge t+2$. It follows that removing the bad column changes
the value of essential differences, but does not change their number.

More generally, $m$ consecutive bad columns occur in $\tilde\Gamma$ if and only if
$\Gamma$ contains $m{+}1$ consecutive simple roots. Here the argument is practically the
same.

(c) \ Assume the column $\varkappa_1=\genfrac{(}{)}{0pt}{}{1}{j_1-1}$
%\left(\!\begin{array}{c} 1 \\ j_1{-}1 \end{array}\!\right)$ 
is bad. Then $j_1=1$ and $i_1=1$, i.e., $\gamma_1=\ap_1$.
If the next-to-right column is good, then $\tilde\Gamma$ contains 
a fragment of the form
\[
  \left(\begin{array}{ccccc} 0 & 1  &     2    & i_2+1 &\dots \\
                                             & 0  & j_2{-}1 & j_3-1  &  \dots
\end{array}\right), 
\]  
where $j_2\ge 3$. Having removed the bad column $\genfrac{(}{)}{0pt}{}{1}{0}$,
%$\left(\begin{array}{c} 1 \\ 1 \end{array}\right)$, 
we gain the essential difference `$2$' in the first row instead of the essential
difference $(j_2{-}1)$ in the second row. However, the total number of essential differences
remains intact. The similar argument applies if there are several consecutive bad columns
including $\varkappa_1$ or if $\varkappa_{k+1}=\genfrac{(}{)}{0pt}{}{i_k+1}{n}$ is bad.
\end{proof}

\noindent In what follows, the function $\eus Y:\AND\to \BN$ is said to be the
{\it OY--invariant}. Here are further properties of $\eus Y$.

\begin{prop}  \label{prom:valuesY}
The minimal (resp. maximal) value of $\eus Y$ is $0$ (resp. $n-1$). 
Each of them is attained on a unique $\fX$-orbit.
Namely, $\eus Y(\Gamma)=0$  if and only if\/ $\Gamma$ lies in the standard $\fX$-orbit;
$\eus Y(\Gamma)=n-1$  if and only if\/ $\Gamma=\{\ap_1,\ap_3,\dots\}$ or
$\{\ap_2,\ap_4,\dots\}$.
\end{prop}\begin{proof}
This is easily verified using Eq.~\eqref{conv_oy}.
\end{proof}

\begin{rmk}
The definition of $\eus Y(\Gamma)$ given in Eq.~\eqref{def-oy} can be repeated 
verbatim for any other root system. However, such a function will not be $\fX$-invariant.
To save $\fX$-invariance, one might attempt to endow summands in  Eq.~\eqref{def-oy}
with certain coefficients. This works in the symplectic case. Namely, one has to put coefficient
`2' in front of $r_\Gamma(\gamma_s)$ if $\gamma_s$ is short. The explanation stems
from the fact that there is an unfolding procedure that takes $\GR{C}{n}$ to $\GR{A}{2n-1}$.
This procedure allows us to identify an antichain (upper ideal) in $\Delta^+(\GR{C}{n})$
with a "self-conjugate" antichain (upper ideal) in $\Delta^+(\GR{A}{2n-1})$, see 
\cite[5.1]{dual} for details. Under this procedure, each short root in $\Delta^+(\GR{C}{n})$
is replaced with two roots in $\Delta^+(\GR{A}{2n-1})$.
Therefore, the modified sum for an antichain  in $\Delta^+(\GR{C}{n})$
actually represents the OY-invariant for the  corresponding "self-conjugate" antichain
in $\Delta^+(\GR{A}{2n-1})$. 
Since $\Delta^+(\GR{C}{n})\simeq \Delta^+(\GR{B}{n})$, the modified formula can also 
be transferred to the $\GR{B}{n}$-setting. But the last isomorphism does not
respect root lengths. Therefore the definition becomes quite unnatural for $\GR{B}{n}$.
\\
Also, it is not clear how to construct an $\fX$-invariant in case of $\GR{D}{4}$.
\end{rmk}

%%%%%%%%%
\subsection{$\fX$-orbits and duality} 
For $\Delta$  of type $\GR{A}{n}$,  we introduced in \cite[\S\,4]{dual} a certain involutory
map (``duality'')   \ $\ast: \AND\to \AND$. 
It has the following properties:

(1) \  $\#\Gamma +\#(\Gamma^*)=n$;

(2)  \ If $\Gamma\subset\Pi$, then $\Gamma^*=\Pi\setminus\Gamma$;

(3) \ $\Delta(i)^*=\Delta(n+2-i)$.

\noindent  Say that $\Gamma^*$ is the {\it dual antichain\/} for $\Gamma$.
Our aim is to establish a relationship between $\fX$ and `$\ast$'. 
To this end, recall the explicit definition of the duality 
$\Gamma\mapsto\Gamma^*$.

Suppose $\Gamma=\{(i_1,j_1),\dots,(i_k,j_k)\}$ as above. In this subsection, we 
represent $\Gamma$ as the usual two-row array:
\[
    \Gamma=\left(\begin{array}{ccc}  i_1 & \dots & i_{k}  \\ 
     j_1 & \dots & j_{k} \end{array}\right) .
\]
Set $\bi=\bi(\Gamma)=(i_1,\dots,i_k)$ and $\bj=\bj(\Gamma)=(j_1,\dots,j_k)$. 
That is, $\Gamma=(\bi,\bj)$ is determined by two strictly 
increasing sequences of equal cardinalities
lying in $[n]:=\{1,\dots,n\}$ such that $\bi\le \bj$ (componentwise).
Then $\Gamma^*=(\bi^*,\bj^*)$ is defined by
\[
      \bi^*:=[n]\setminus \bj \ \text{ and } \ \bj^*:=[n]\setminus \bi \ .
\]
It is not hard to verify that $\Gamma^*$ is an antichain, see \cite[Theorem\,4.2]{dual}
(Our notation for the roots of $\slno$ is slightly different from that in \cite{dual}, therefore
the definition of $\Gamma^*$ has become a bit simpler.)

\begin{thm}
For any $\Gamma\in \AND$, we have $\fX(\Gamma)^*=\fX^{-1}(\Gamma^*)$.
\end{thm}\begin{proof}
We prove that the $\bi$- and $\bj$-sequences for $\fX(\Gamma)^*$ and
$\fX^{-1}(\Gamma^*)$ coincide.

\noindent Below, we use the description of $\fX$ given in the proof of Theorem~\ref{thm:osya}.
We have
\[
    i\in \bi(\fX(\Gamma)^*) \Leftrightarrow i \not\in \bj(\fX(\Gamma)) \Leftrightarrow 
\begin{cases}  i+1\not\in \bj(\Gamma) \text{ or }\\
                       i+1\in \bj(\Gamma) \text{ and } \ \ap_i,\ap_{i+1}\in \Gamma .
\end{cases}                   
\]
The last possibility means that $i+1$ occurs in the $\bj$-sequence of $\Gamma$ and hence
$i$ occurs in the $\bj$-sequence of $\tilde\Gamma$; however, it occurs in a bad column and
%$\tilde\Gamma$ contains bad columns, and the number $i+1$ 
therefore disappears after removing the bad columns.

On the other hand,  consider $\bi(\fX^{-1}(\Gamma^*))$. To this end, one needs an
explicit description of $\fX^{-1}$ in terms of two-row arrays.
As the description of $\fX$ includes deletion of some columns, the description of 
$\fX^{-1}(\Gamma^*)$ should include a creation of columns. 
More precisely, for 
\[
\Gamma^*=\left(\begin{array}{ccc}  i_1^* & \dots & i_{n-k}^*  \\ 
     j_1^* & \dots & j_{n-k}^* \end{array}\right) ,
\]
we perform the following.
First, if $i_1^*\ge 2$, then we put columns 
$\left(\begin{array}{ccc}  1 & \dots & i_{1}^*{-}1  \\ 
     1 & \dots & i_{1}^*{-}1 \end{array}\right)$ at the beginning.
Then each pair of consecutive columns of $\Gamma^*$ is transformed  as follows:
\[
    \left(\begin{array}{cc}  i_s^* &  i_{s+1}^*  \\ 
     j_s^*  & j_{s+1}^*  \end{array}\right)\mapsto
    \begin{cases} \left(\begin{array}{c}  i_{s+1}^*{-}1 \\ j_s^*{+}1 \end{array}\right)
      &  \text{if } \ i_{s+1}^*{-}1\le j_s^*{+}1, \\
    \left(\begin{array}{ccc}   j_s^*{+}1& \dots &  i_{s+1}^*{-}1  \\ 
     j_s^*{+}1& \dots &  i_{s+1}^*{-}1 \end{array}\right) &  \text{if } \ i_{s+1}^*{-}1> j_s^*{+}1.
     \end{cases}
\]
Finally,  if $ j_{n-k}^*< n$, then we put columns 
$\left(\begin{array}{ccc}   j_{n-k}^*{+}1 & \dots & n  \\ 
     j_{n-k}^*{+}1 & \dots & n \end{array}\right)$ at the end.
The resulting two-row array represents $\fX^{-1}(\Gamma^*)$. From this description, it follows
that

\begin{multline*}
    i\in \bi(\fX^{-1}(\Gamma^*)) \Longleftrightarrow 
    \begin{cases}  i+1\in \bi(\Gamma^*) \text{ or } \\
                      %\text{ or } \   i+1\not\in \bi(\Gamma), \text{ but } \ i< i^*_1-1 \\
                i+1\not\in \bi(\Gamma^*) \text{ but} \ 
                \begin{cases}  i \le i^*_1-1  \text{ \ or} \\       
                                        i\ge j^*_{n-k}+1 \text{ or} \\         
                 j^*_s{+}1\le i \le  i^*_{s+1}-1 \text{ for some } s\in [n{-}k{-}1] 
                                      \end{cases}
\end{cases}   \\
\Longleftrightarrow   \begin{cases} i+1\not\in \bj(\Gamma) \text{ or } \\ 
                       i+1\in \bj(\Gamma) \text{ and } \ \ap_i,\ap_{i+1} \in\Gamma.                
\end{cases}   
\end{multline*} 
Thus, we have proved that $\bi(\fX(\Gamma)^*)=\bi(\fX^{-1}(\Gamma^*))$. The argument for 
$\bj$-sequences is similar.
\end{proof}

There is also a connection between the duality and OY--invariant:

\begin{prop}
$\eus Y(\Gamma)=\eus Y(\Gamma^*)$.
\end{prop}\begin{proof}
As above, we think of $\Gamma$ as union of sequences $I=(i_1,\dots,i_k)$ and 
$J=(j_1,\dots,j_k)$.
Using Eq.~\eqref{conv_oy}, we write
\[
  \eus Y(\Gamma)=r_\bullet(I)+ r^\bullet(J) ,
\]
where $r_\bullet(I)=\sum_{s=1}^k\chi(i_s{-}i_{s-1})$ and 
  $r^\bullet(J)=\sum_{s=1}^k \chi(j_{s+1}{-}j_s)$. Recall that $i_0=0$ and $j_{k+1}=n+1$.
  Then the assertion will follow from
the definition of $\Gamma^*$ and the equalities
$r_\bullet(I)=r^\bullet([n]\setminus I)$ and $r^\bullet(J)=r_\bullet([n]\setminus J)$.
Clearly, it suffices to prove one of them.

Let us say that $C_i=\{c_i,c_i+1,\dots, c+{m_i}\}$ is a {\it connected component\/} of 
$I\cup\{0\}$, if  $C_i\subset I\cup\{0\}$ and $c_i-1, c+m_i+1\not\in I\cup\{0\}$.
One similarly defines the connected components of $J\cup\{n{+}1\}$.
Since the consecutive differences inside a connected component are unessential,
we obtain
\begin{gather*}
 r_\bullet(I)=\bigl(\text{the number of connected components of $I\cup\{0\}$}\bigr)-1 , \\
 r^\bullet(J)=\bigl(\text{the number of connected components of $J\cup\{n{+}1\}$}\bigr)-1 .
\end{gather*}
Now, the equality $r_\bullet(I)=r^\bullet([n]\setminus I)$ can be proved using a simple
verification. One has to consider four cases depending on whether $1$ and $n$ belong 
to $I$. As a sample, we consider one case.

Assume $1,n\not\in I$. Then $\{0\}$  is a connected component of $I\cup\{0\}$. 
If $I$ itself has $m$ connected components, then the total number of components is $m+1$.
Hence $r_\bullet(I)=m$. On the other hand, the assumption shows that
$[n]\setminus I$ has $m+1$ connected components. Furthermore, 
$n\in ([n]\setminus I)$. Therefore
$\{n+1\}$ does not form a connected component. Thus, $([n]\setminus I)\cup\{n{+}1\}$ 
still has $m+1$ components, and $r^\bullet ([n]\setminus I)=m$.
\end{proof}

%%%%%%%%%%%%%%
\appendix
\section{Computations for $\GR{F}{4}$}
%%%%%%%%%%%%%%

\noindent  We use the numbering of simple roots from \cite{vo}.
The positive root $\beta=\sum_{i=1}^4 n_i\ap_i$ is denoted by
$(n_1n_2n_3n_4)$. For instance, $\theta=(2432)$ and $\theta_s=(2321)$.

I. \ $\#\AND=105$ and $h=12$. There are eleven  $\fX$-orbits:
eight orbits of cardinality 12 and orbits of cardinality 2,3, and 4.
We indicate representatives and cardinalities  for all orbits:

 \ \{1000\}   -- 12;  \ \{0100\}  -- 12; \ \{0010\} -- 12; \ \{0001\} -- 12;
 \ \{0011\} -- 12;
 
\ \{1100\} -- 12;
\ \{1111\} -- 12; 
\ \{2432\} -- 12  (the standard orbit);
\ \{1000, 0010\}  -- 2;

\ \{0110\}  -- 3;
\ \{0001, 1110\} -- 4.

\vskip1ex

II. \ $\#\ANDo=66$ and $h-1=11$. The notation $\Gamma\leadsto\Gamma'$ means
$\Gamma'=\fX_0(\Gamma)$.
The $\fX_0$-orbits are:

\noindent
1) \ The standard one: $\Delta(11)=\{2432\} \leadsto \{2431\}\leadsto \dots \leadsto \Delta(2)\leadsto 
\varnothing\leadsto\Delta(11)$;
\\[.7ex]
2) \ $\{1321\}\!\leadsto\! \{2221\} \!\leadsto\! \{1321,2211\} \!\leadsto\! \{1221, 2210\} \!\leadsto\! \{0221,1211\} \leadsto
\{0211, 1111, 2210\} \!\leadsto\! \{0111, 1210\} \!\leadsto\! \{0011,0210, 1110\} \!\leadsto\! \{0110, 1100\}
\!\leadsto\! \{0011\} \!\leadsto\! \{2210\} \!\leadsto\!  \{1321\}$;
\\[.7ex]
3) $\{1221\}\!\leadsto\! \{0221, 2211\}\!\leadsto\! \{1211, 2210\}\!\leadsto\! \{0221, 1111, 1210\}\!\leadsto\! \{0211, 1110\}
\!\leadsto\!  \{0111, 0210, 1100\}\!\leadsto\! \{0011, 0110\}\!\leadsto\! \{1100\}\leadsto\{0221\}\!\leadsto\! \{2211\}
\!\leadsto\! \{1321, 2210\}\!\leadsto\! \{1221\}$;
\\[.7ex]
4) \ $\{1211\}\!\leadsto\! \{0221, 1111, 2210\}\!\leadsto\! \{0211, 1210\}\!\leadsto\! \{1111,0210\}\!\leadsto\! \{0111,1110\}\leadsto
\{0011, 0210, 1100\}\!\leadsto\! \{0110\}\!\leadsto\! \{0011, 1100\}\!\leadsto\! \{0210\}\!\leadsto\! \{1111\}\!\leadsto\! \{0221, 2210\}\leadsto
\{1211\}$;
\\[.7ex]
5) $\{1210\}\!\leadsto\! \{0221, 1111\}\!\leadsto\! \{0211, 2210\}\!\leadsto\! \{1111, 1210\}\!\leadsto\! \{0221, 1110\}\!\leadsto\! \{0211, 1100\}
\!\leadsto\! \{0111, 0210\}\!\leadsto\! \{0011, 1110\}\!\leadsto\! \{0210, 1100\}\!\leadsto\! \{0111\}\!\leadsto\! \{0011, 2210\}\!\leadsto\! \{1210\}
$;
\\[.7ex]
6) $\{1110\}\!\leadsto\! \{0221, 1100\}\!\leadsto\! \{0211\}\!\leadsto\! \{1111, 2210\}\!\leadsto\! \{0221, 1210\}\!\leadsto\! \{0211, 1111\}\leadsto
\{0111, 2210\}\!\leadsto\! \{0011, 1210\}\!\leadsto\! \{0210, 1110\}\!\leadsto\! \{0111, 1100\}\!\leadsto\! \{0011, 0210\}\!\leadsto\! \{1110\}
$.

Each orbit consists of 11 antichains.

III. \ $\#\ANDs=21$ and $h^*=9$. The $\fX_s$-orbits are:

\noindent
1) \  standard: $\Delta_s(8)=\{2321\}\!\leadsto\! \{1321\} \!\leadsto\!\dots \leadsto\! \Delta_s(1)=\{1000, 0100\} \leadsto \varnothing \!\leadsto\! \Delta_s(8)$;%\{\theta_s\}$;
\\[.7ex]
2) \ $\{0100\}\leadsto \{1000\}\leadsto \{0111\}\leadsto \{1210\}\leadsto \{1111\}\leadsto \{0111, 1210\}\leadsto \{1110\}\leadsto 
\{0111, 1100\}\leadsto \{0110, 1000\}\leadsto \{0100\}$;
\\[.7ex]
3) \ $\{1100\}\leadsto \{0111, 1000\}\leadsto \{0110\}\leadsto \{1100\}$.

\end{document}